\documentclass{article}
\usepackage{amsthm, amsfonts, amssymb, latexsym, amsmath}
\usepackage{color}

\usepackage[utf8]{inputenc}
\newcommand{\F}{\mathbb{F}}
\newcommand{\Z}{\mathbb{Z}}

\newtheorem{theorem}{Theorem}
\newtheorem{lemma}{Lemma}

\title{A Note on the Cross-Correlation of Costas~Permutations}
\author{Domingo Gomez-Perez$^1$ and Arne Winterhof$^2$\\
$^1$ Facultad de Ciencias Universidad de Cantabria,\\ 
Santander, Spain\\
E-mail: domingo.gomez@unican.es\\
$^2$ Johann Radon Institute for\\ Computational and Applied Mathematics,\\
Austrian Academy of Sciences, Altenberger Str. 69,\\ 4040 Linz, Austria\\
E-mail: arne.winterhof@oeaw.ac.at
}
\date{}

\begin{document}

\maketitle

\begin{abstract}
    We build on the work of Drakakis et al.\ (2011) on the maximal cross-correlation of the families
    of Welch and Golomb Costas permutations. In particular, we settle some of their conjectures.
    More precisely, we prove two results. 
    
    First, for a prime $p\ge 5$, the maximal cross-correlation of the family of the $\varphi(p-1)$ different Welch Costas permutations of $\{1,\ldots,p-1\}$ is $(p-1)/t$, where $t$ is the smallest prime divisor of $(p-1)/2$ if $p$ is not a safe prime and at most $1+p^{1/2}$ otherwise. Here $\varphi$ denotes Euler's totient function and a prime $p$ is a safe prime if $(p-1)/2$ is also prime.
    
    Second, for a prime power $q\ge 4$ the maximal cross-correlation of a subfamily of Golomb Costas permutations of $\{1,\ldots,q-2\}$ is $(q-1)/t-1$ if~$t$ is the smallest prime divisor of $(q-1)/2$ if $q$ is odd and of $q-1$ if~$q$ is even provided that $(q-1)/2$ and $q-1$ are not prime, and at most $1+q^{1/2}$ otherwise.
        Note that we consider a smaller family than Drakakis et al. Our family is of size $\varphi(q-1)$ whereas there are $\varphi(q-1)^2$ different Golomb Costas permutations. The maximal cross-correlation of the larger family given in the tables of Drakakis et al.\ is larger than our bound (for the smaller family) for some $q$.
        \end{abstract}



\section{Introduction}

For a positive integer $n$, let $\pi$ be a permutation of $\{1,\ldots,n\}$ satisfying
$$\pi(i+k)-\pi(i)\ne \pi(j+k)-\pi(j)$$ 
for any integers $1\le k\le n-2$ and $1\le i<j\le n-k$.
Such a permutation is called a {\em Costas permutation} of $\{1,\ldots,n\}$
and the corresponding $(n\times n)$-permutation matrix $A=(a_{ij})_{i,j=1}^n$ defined by
$$a_{ij}=1 \mbox{ if and only if }\pi(i)=j$$
is called a {\em Costas array} of size $n$.
These objects are crucial in some problems arising from radar and sonar, see for example \cite[Section~7.6]{ju} and \cite{gota}.

The {\em cross-correlation $C_{f_1,f_2}(u,v)$} between two mappings 
$$f_1,f_2:\{1,\ldots,n\}\rightarrow \{1,\ldots,n\}$$
at $(u,v)\in \Z^2$, $1-n\le u,v\le n-1$, is the number of solutions 
$$x\in \{\max\{1,1-u\},\ldots,\min\{n,n-u\}\}$$ 
of the equation
\begin{equation}\label{cor} f_1(x)+v=f_2(x+u).
\end{equation}
For a family ${\cal F}$ of Costas permutations of $\{1,\ldots,n\}$, the {\em maximal cross-correlation} $C({\cal F})$ is
$$C({\cal F})=\max_{u,v} \max_{\genfrac{}{}{0pt}{}{f_1,f_2\in {\cal F}}{f_1\ne f_2}}C_{f_1,f_2}(u,v).$$
Studying the maximal cross-correlation of a family of Costas permutations is not only a very interesting mathematical problem, since families with small maximal cross-correlation are of high practical importance, see \cite{drakakis2011maximal} and references therein.

In this note, we study the maximal cross-correlation of two families of Costas permutations, the
family of Welch Costas permutations and a subfamily of Golomb Costas permutations defined below. In particular, we will address some open problems from \cite{drakakis2011maximal}.

{\em Welch's construction} of Costas permutations is defined as follows, see \cite{gota,ju}. For a prime $p>2$, let $g$ be a primitive root modulo $p$ and $\pi_g$ the permutation of $\{1,\ldots,p-1\}$ defined by
$$\pi_g(i) \equiv g^i \bmod p.$$
Then, for $p\ge 5$, the family ${\cal W}_p$ of Welch Costas permutations of $\{1,\ldots,p-1\}$ is
$${\cal W}_p=\{\pi_g : g \mbox{ primitive root modulo $p$}\},$$
so that, 
$|{\cal W}_p|=\varphi(p-1),$
where $\varphi$ is {\em Euler's totient function}.

A prime $p$ is a {\em safe prime} if $(p-1)/2$ is also 
a prime, called {\em Sophie Germain prime}.
Therefore,
$$|{\cal W}_p|=\frac{p-3}{2}\quad \mbox{if $p\ge 7$ is a safe prime}$$
and $|{\cal W}_5|=2$.

In this note we prove the following result on $C({\cal W}_p)$.

\begin{theorem}
\label{main}
For a prime $p\ge 5$, let $t$ be the smallest prime divisor of \hbox{$(p-1)/2$}. Then, the maximal cross-correlation $C({\cal W}_p)$ of the family of Welch Costas permutations ${\cal W}_p$ of $\{1,\ldots,p-1\}$ satisfies
$$C({\cal W}_p)\left\{\begin{array}{ll} \le 1+\lfloor(1-2/(p-1))p^{1/2}\rfloor & \mbox{if $p$ is a safe prime},\\
= (p-1)/t & \mbox{otherwise}.\end{array} \right.$$
\end{theorem}
Note that we can substitute each $\pi_g(i)$ by a shift $\pi_g(i+c_g)$ and get the same result.
However, ${\cal W}_p$ must not contain two shifts for the same primitive element~$g$.
In particular, for non-safe primes, Theorem~\ref{main} settles the first conjecture in
Drakakis et al.\ \cite[Conjecture 3]{drakakis2011maximal}.
We prove Theorem~\ref{main} in Section~\ref{secW}.

%

{\em Golomb's construction} of Costas permutations is the following, see \cite{comu,Golomb1984,ju,moso}. For a prime power $q>2$ and primitive elements $g_1$ and $g_2$ of the finite field $\F_q$, let $\pi_{g_1,g_2}$ be the permutation of $\{1,\ldots,q-2\}$ defined by 
$$\pi_{g_1,g_2}(i)=h \mbox{ if and only if } g_1^i+g_2^h=1.$$
For $q\ge 4$ and fixed $g_2$ we study the subfamily ${\cal G}_q$ of the family of Golomb Costas permutations of $\{1,\ldots,q-2\}$ defined by
$${\cal G}_q=\{\pi_{g_1,g_2} : g_1 \mbox{ primitive element of $\F_q$}\}.$$
Then we have $|{\cal G}|_q=\varphi(q-1).$
In Section~\ref{golomb}, we prove the following result on $C({\cal G}_q)$.
\begin{theorem}\label{main2}
 For a prime power $q\ge 4$, let $t$ be the smallest prime divisor of \hbox{$(q-1)/2$} if $q$ is odd and
 of $q-1$ if $q$ is even. Then, the maximal cross-correlation $C({\cal G}_q)$ of the family of Golomb Costas permutations ${\cal G}_q$ of $\{1,\ldots,q-2\}$ satisfies
 $$C({\cal G}_q)\left\{\begin{array}{ll}\le 1+\lfloor (1-2/(q-1))q^{1/2}\rfloor & \mbox{if $q$ is odd and }t=(q-1)/2,\\
 \le \lfloor(1-1/(q-1))(1+q^{1/2})\rfloor & \mbox{if $q$ is even and $t=q-1$},\\
 = (q-1)/t-1 & \mbox{otherwise}.\end{array}\right.$$
\end{theorem}

Besides $C({\cal G}_q)$, it is interesting to study the cross-correlation $C({\cal L}_q)$ 
of the larger set ${\cal L}_q$ of all Golomb Costas permutations
$${\cal L}_q=\{\pi_{g_1,g_2}: g_1,g_2 \mbox{ primitive elements of }\F_q\}$$
of size
$$|{\cal L}_q|=\varphi(q-1)^2.$$
The tables of \cite{drakakis2011maximal} show that $C({\cal L}_q)$ is larger than $C({\cal G}_q)$ for some small values of $q$.
For example, for $q=59$, we have $C({\cal L}_{59})=12$ but $C({\cal G}_{59})\le 8$. 
However, for all prime values of $q$ with $61\le q\le 271$ and all strict prime powers $25\le q\le 343$, the bound of Theorem~\ref{main2} is also valid for $C({\cal L}_q)$.
It remains an open problem to prove the conjecture that this bound holds for $C({\cal L}_q)$ up to a few exceptions of $q$ with~$q\le 59$.

\section{Proof of Theorem~\ref{main}}
\label{secW}

By \cite[Theorem~1]{drakakis2011maximal}, we have
$$\max_{u\in \Z} \max_{\genfrac{}{}{0pt}{}{f_1,f_2\in {\cal W}}{f_1\ne f_2}}C_{f_1,f_2}(u,0)=\frac{p-1}{t}.$$
Since $t\le \sqrt{(p-1)/2}$ if $p$ is not a safe prime, it remains to prove the following lemma, from which Theorem~\ref{main} follows immediately
after verifying
$$\frac{p-1}{t}\ge \sqrt{2(p-1)}\ge 1+p^{1/2}\quad \mbox{for }p\ge 11$$
and that $5$ and $7$ are both safe primes.

\begin{lemma}\label{Welch}
For any prime $p\ge 5$ we have
$$\max_u \max_{v\ne 0}\max_{\genfrac{}{}{0pt}{}{f_1,f_2\in {\cal W}_p}{f_1\ne f_2}}C_{f_1,f_2}(u,v)\le 1+\left\lfloor\left(1-\frac{2}{p-1}\right)p^{1/2}\right\rfloor.$$
\end{lemma}
Proof.
The maximum in the statement can be bounded by the maximal number $N$ of solutions 
$x\in \F_p^*$ of any equation of the form
\begin{equation}\label{eq} ax^r\equiv x+v \bmod p,\quad \mbox{with } av\not\equiv 0\bmod p,\quad
\gcd(r,p-1)=1,~1<r<p-1,
\end{equation}
since, if $g$ is a fixed primitive root modulo $p$, all other primitive roots modulo~$p$
are of the form $g^r$ with $\gcd(r,p-1)=1$.
For fixed $a$ and $v$ with $av\not\equiv 0\bmod p$, the number of solutions of $(\ref{eq})$
is
$$\frac{1}{p-1}\sum_{\chi} \sum_{x\in \F_p^*\setminus\{-v\}} \chi(ax^r)\overline{\chi}(x+v)$$
by the orthogonality relations
$$\frac{1}{p-1}\sum_{\chi}\chi(x)\overline{\chi}(y)=\left\{\begin{array}{cc} 1, & x=y,\\ 0, & x\ne y,\end{array}\right\} \quad\mbox{for all } x,y \in \F_p^*,$$
where the sum runs through all multiplicative characters $\chi$ of $\F_p$.

The contribution of the trivial character $\chi_0$ is $(p-2)/(p-1)$ and that of the quadratic character $\eta$ is
$-\eta(a)/(p-1)$ by \cite[Lemma 7.3.7]{ju}. 
Thus, 
\begin{eqnarray*} 
N&\le& 1+\frac{p-3}{p-1}\max_{v\in \F_p^*}\max_{\chi\not\in\{\chi_0,\eta\}}\left|\sum_{x\in \F_p^*}\chi(x^r(x+v)^{p-2})\right|\\
&\le& 1+\left(1-\frac{2}{p-1}\right)p^{1/2}
\end{eqnarray*}
by the Weil bound, see for example \cite[Theorem 5.41]{LN}.~\hfill $\Box$

\section{Proof of Theorem~\ref{main2}}
\label{golomb}
For $u=v=0$, we have, by \cite[Theorem~3]{drakakis2011maximal}
$$ \max_{\genfrac{}{}{0pt}{}{f_1,f_2\in {\cal G}_q}{f_1\ne f_2}}C_{f_1,f_2}(0,0)=\frac{q-1}{t}-1,
$$
where $t$ is the smallest prime divisor of $(q-1)/2$ if $q$ is odd and of $q-1$ if $q$ is even.

Next we prove an upper bound for $v=0$ and arbitrary $u$.
\begin{lemma}\label{gol0}
 We have
 $$\max_{u} \max_{\genfrac{}{}{0pt}{}{f_1,f_2\in {\cal G}_q}{f_1\ne f_2}}C_{f_1,f_2}(u,0)\le
  \frac{q-1}{t}-1$$
  if $t\not\in\{(q-1)/2,q-1\}$
  and
  $$\max_{u} \max_{\genfrac{}{}{0pt}{}{f_1,f_2\in {\cal G}_q}{f_1\ne f_2}}C_{f_1,f_2}(u,0)\le 2$$
  otherwise.
\end{lemma}
Proof.
Since 
$$C_{f_1,f_2}(-u,0)=C_{f_2,f_1}(u,0)$$ 
we may assume $u\ge 1$.
Let $f_1$ and $f_2$ be defined by $f_1(x)=h$ if and only if $g_1^x+g_2^h=1$ and 
$f_2(x)=h$ if and only if $g_1^{xr}+g_2^h=1$, respectively, for some integer $r$ with $\gcd(r,q-1)=1$ and
$1<r<q-1$. Then, the number of solutions~$x$ of $(\ref{cor})$ with $v=0$ (and $n=q-2$) is the number of integers $x$ in the range $1\le x\le q-2-u$ such that
$$ g_1^x=g_1^{(x+u)r},$$
that is, $x$ satisfies
$$(r-1)x\equiv -ur \bmod (q-1).$$
Put $d=\gcd(r-1,q-1)$ and let $a$ be the inverse of $(r-1)/d$ modulo $(q-1)/d$.
There is no solution if $d$ does not divide~$u$. Otherwise, the solutions
are those $x$ with 
\begin{equation}\label{equal}x\equiv -a(u/d)r \bmod (q-1)/d.
\end{equation}
We have at most $d$ such solutions $x$ with $1\le x\le q-2$. 
Obviously, we have $d\le (q-1)/t$.
The result follows immediately if $d< (q-1)/t$. It remains to study the case
$d=(q-1)/t$. 
Then either 
$$u\ge d=(q-1)/t\ge (q-1)^{1/2}\ge t$$ or 
$$\mbox{for $q$ odd}, \quad t=(q-1)/2 \mbox{ and }d=2 $$
and 
$$\mbox{for $q$ even}, \quad t=q-1 \mbox{ and }d=1.$$ 
In the fist case, the solutions $x$ of~$(\ref{equal})$ are of the form $x=x_0+kt$ with $1\le x_0\le t$ and
$0\le k\le d-1$. However, $k=d-1$ is not possible since 
$$x_0+(d-1)t> q-1-t>q-2-(q-1)^{1/2}\ge q-2-d\ge q-2-u$$
and there are at most $d-1=(q-1)/t-1$ solutions.
In the remaining cases we have at most $2$ solutions.
~\hfill $\Box$\\

For $v\ne 0$, analogously to Lemma~\ref{Welch}, we get the following bound.
\begin{lemma}\label{gol1}
For odd $q$, we have
$$\max_u \max_{v\ne 0} \max_{\genfrac{}{}{0pt}{}{f_1,f_2\in {\cal G}_q}{f_1\ne f_2}}(u,v)\le 1+\left\lfloor\left(1-\frac{2}{q-1}\right)q^{1/2}\right\rfloor$$
and, for even $q$,
$$\max_u \max_{v\ne 0}\max_{\genfrac{}{}{0pt}{}{f_1,f_2\in {\cal G}_q}{f_1\ne f_2}}(u,v)\le \left\lfloor\left(1-\frac{1}{q-1}\right)\left(1+q^{1/2}\right)\right\rfloor.$$
\end{lemma}
Proof. 
Again, let $f_1(x)=h$ whenever $g_1^x+g_2^h=1$ and $f_2(x)=h$ whenever $g_1^{rx}+g_2^h=1$
for some $r$ with $\gcd(r,q-1)=1$ and $1<r<q-1$. Then, $(\ref{cor})$ implies
$$g_2^v(1-g_1^x)=1-g_1^{r(x+u)}.$$
Substituting $y=1-g_1^x$, $a=g_2^v$ and $b=g_1^{ru}$, we get
$$ay=1-b(1-y)^r.$$
Note that $a\ne 1$ since $v\not= 0$ and $y\not\in\{0,1\}$ since $1\le x\le q-2$.
Hence, we have to estimate the number $N$ of solutions of equations of the form
$$b(1-y)^r=1-ay,\quad y\in \F_q^*\setminus\{1\},$$
for any $a\in \F_q^*\setminus\{1\}$ and $b\in \F_q^*.$
We can represent $N$ by character sums
$$N=\frac{1}{q-1} \sum_{\chi} \sum_{y\in \F_q^*\setminus\{1\}} \chi(b) \chi((1-y)^r(1-ay)^{q-2}).$$
The contribution of the trivial character is $(q-2)/(q-1)$ and, for odd $q$, that of
the quadratic character is $-\chi(b)/(q-1)$. For the remaining characters, the absolute value of the inner sum is at most $q^{1/2}$.
Collecting these facts, the result follows.~\hfill $\Box$\\

Theorem~\ref{main2} is proved by combining Lemmas~\ref{gol0} and \ref{gol1}, 
after verifying that $t<(q-1)/2$ implies the following results.

For any odd $q$, we have $t\le \sqrt{\frac{q-1}{2}}$, and thus 
$$1+q^{1/2}\le \frac{q-1}{t}-1 \quad \mbox{for any odd $q\ge 27$}.$$
For the remaining odd $q$ with $q\le 25$, the following refinement holds, 
$$1+\left\lfloor\left(1-\frac{2}{q-1}\right)q^{1/2}\right\rfloor\le \frac{q-1}{t}-1.$$

For even $q\ge 4$, by Mih\u ailescu’s Theorem (former Catalan conjecture), $q-1$ is not a perfect square and thus $q-1\ge t(t+2)$, that is, $t\le -1+q^{1/2}$.

If $q=2^r$ with an odd $r$, then 
$$1+\lfloor q^{1/2}\rfloor < 1+ q^{1/2}\le \frac{q-1}{t}.$$
If $q=2^r$ with an even $r$, then $t=3$ and
$$1+q^{1/2}\le \frac{q-1}{3}-1 \quad \mbox{for $q\ge 64$}.$$
In the remaining case, that is $q=16$, the more precise bound of Lemma 3 equals $(q-1)/t-1=4$.~\hfill $\Box$

\section*{Acknowledgments}

The first author is partially supported 
by the Consejer\'ia de Universidades e Investigaci\'on, Medio Ambiente y Politica Social del Gobierno de Cantabria, reference number Vp34.

The second author is partially supported by the Austrian Science Fund FWF Project P 30405-N32.

\end{document}